\documentclass[12pt]{article}

\usepackage[english]{babel}
\usepackage{amssymb,amsmath}

\newtheorem{theorem}{Theorem }[section]

\newtheorem{corollary}[theorem]{Corollary}

\newcommand{\proof}{\noindent\textbf{Proof. }}
\newcommand{\qed}{\hspace*{\fill}$\Box$}

\providecommand{\keywords}[1]{\textbf{Keywords---} #1}
\providecommand{\AMS}[1]{\textbf{AMS classification---} #1}

\def\srg{\mathop{\mathrm {srg}}\nolimits}

\newcommand{\cD}{\mathcal{D}}

\newcommand{\cG}{\mathcal{G}}

\title{Automorphisms of Strongly Regular Graphs}

\author{S. De Winter, E. Kamischke, Z. Wang}

\begin{document}

\maketitle

\abstract{In this article we generalize a theorem of Benson for generalized quadrangles to strongly regular graphs and directed strongly regular graphs. The main result provides numerical restrictions on the number of fixed vertices and the number of vertices  mapped to adjacent vertices under an automorphism. It is explained how these results can be used when studying partial difference sets in Abelian groups and projective two-weight sets. The underlying ideas are linear algebraic in nature.}
\medskip


\AMS{05C50, 05E30}

\section{Introduction}

In 1970 Benson \cite{Benson} provided a congruence that relates the parameters of a finite generalized quadrangle  to the number of fixed points and the number of points mapped to collinear points under an automorphism. This result has turned out to be useful at various occasions in the theory of generalized quadrangles, see for example \cite{Payne-Thas, SDW-KT}. In 2006 the first author \cite{SDW} generalized this theorem to partial geometries and used it to obtain a characterization of the so-called Van Lint - Schrijver partial geometry.  In 2010 Temmermans in her dissertation (see also \cite{HVM-T-Thas}) provided further generalizations for other geometries, including partial quadrangles, and used these generalizations to study polarities of the geometries under investigation. Though interesting results were obtained two observations stood out at this point that formed the starting point of this paper. On the one hand the Benson type results for generalized quadrangles and partial geometries were particularly elegant, whereas the result for partial quadrangles for example  did not provide an elegant congruence, but rather a more complicated equation. This is due to the fact that a key matrix used in the proof has an eigenvalue equal to zero in the former cases, but not in the latter case. On the other hand, given the abundance  of Benson type theorems for geometries whose collinearity graph is strongly regular one would expect the existence of a unifying theorem for strongly regular graphs. This paper will provide such unifying theorem while at the same time overcoming the problem that arose when the aforementioned matrix has no eigenvalue equal to zero. The aim of the article is not to provide a multitude of new applications, but rather to provide the theory, and point out some  areas in which  these results can be used through some basic examples.

\section{Main Equalities and Congruences}

\subsection{Strongly regular graphs}

We assume the reader to be familiar with strongly regular graphs. For more on these graphs see for example \cite{Godsil}.
Let $\cG$ be a strongly regular graph $\srg(v,k,\lambda,\mu)$, and let $A$ be its adjacency matrix. Then the $v\times v$-matrix $A$ has eigenvalues 
$$\nu_1:=k,$$ 
$$\nu_2:=\frac{1}{2}(\lambda-\mu+\sqrt{\Delta)}),$$ 
$$\nu_3:=\frac{1}{2}(\lambda-\mu-\sqrt{\Delta}),$$ 
with respective multiplicities $$m_1:=1,$$ 
$$m_2:=\frac{1}{2}\left(v-1-\frac{2k+(v-1)(\lambda-\mu)}{\sqrt{\Delta}}\right) $$ 
and $$m_3:=\frac{1}{2}\left(v-1+\frac{2k+(v-1)(\lambda-\mu)}{\sqrt{\Delta}}\right),$$
where $\Delta=(\lambda-\mu)^2+4(k-\mu)=(\nu_2-\nu_3)^2$.
Furthermore these eigenvalues are integers if $2k+(v-1)(\lambda-\mu)\neq0$, that is, if $\cG$ is not a conference graph. In the case where $\cG$ is a conference graph, the eigenvalues are still  integers provided $v$ is a perfect square. We will call these eigenvalues and multiplicities also the eigenvalues and multiplicities of the strongly regular graph.
\medskip

Let $\phi$ be an automorphism of order $n$ of the graph $\cG$. Then $\phi$ corresponds to a $v\times v$-permutation matrix $P$ with the property that $PAP^T=A$. As $P^{-1}=P^T$ we have $PA=AP$. Also $P^n=I$, where $I$ is the $v\times v$-identity matrix. Let $f_\phi$ denote the number of vertices fixed by $\phi$, and let $g_\phi$ be the number of vertices mapped to adjacent vertices by $\phi$ (this excludes fixed vertices). We will write $f$ and $g$ instead of $f_\phi$ and $g_\phi$ if the automorphism $\phi$ is understood. Finally we will also need the following result. If $P$ is a permutation matrix of order $n$, then the eigenvalues of $P$ are the $n$th roots of unity, where the multiplicity of a given primitive $d$th root of unity equals $\sum_{d\mid l} c_l$, where $c_l$ is the number of cycles of length $l$ in the disjoint cycle decomposition of the permutation $P$.

\begin{theorem}\label{MAIN} Let $\cG$ be a strongly regular graph $\srg(v,k,\lambda,\mu)$ whose adjacency matrix $A$ has integer eigenvalues $k$, $\nu_2$ and $\nu_3$. Let $\phi$ be an automorphism of order $n$ of $\cG$, and let $\mu()$ be the M\"obius function. Then for every integer $r$ and all positive divisors $d$ of $n$, there are non-negative integers  $a_d$ and $b_d$ such that 
\begin{equation}k-r +  \sum_{d\mid n} a_d\mu(d)(\nu_2-r) +  \sum_{d\mid n} b_d\mu(d)(\nu_3-r)=-rf+g,\end{equation}
where $f$ is the number of fixed vertices of $\phi$ and $g$ is the number of vertices that are adjacent to their image under $\phi$. Furthermore $a_1+b_1=c-1$, where $c$ is the number of cycles in the disjoint cycle decomposition of $\phi$, and $a_d+b_d=\sum_{d\mid l} c_l$, $d\neq1$, where $c_l$ is the number of cycles of length $l$ of $\phi$.
As a consequence  the following congruence holds:
\begin{equation}k-\nu_3\equiv -\nu_3f+g \pmod{ \sqrt{\Delta}}.\end{equation}
\end{theorem}
\proof Let $M$ be the matrix $M=A-rI$. Then obviously $M$ has integer eigenvalues $\tau_1=k-r$, $\tau_2=\nu_2-r$ and $\tau_3=\nu_3-r$, with respective multiplicities $m_1$, $m_2$ and $m_3$. If $P$ is the permutation matrix corresponding to $\phi$, then clearly $PM=MP$, and hence $(PM)^n=P^nM^n=M^n$. It follows that the eigenvalues of $PM$ are the eigenvalues of $M$ multiplied with appropriate $n$th roots of unity.  As the sum of the elements in each row of $M$ equals $k-r$, the same will hold for $PM$, and hence $k-r$ is an eigenvalue of $PM$.  This eigenvalue clearly has multiplicity $m_1=1$.   Now let $d$ be a positive divisor of $n$, and let $\xi_d$ be a primitive $d$th root of unity. As the eigenvalues of $M$ are integers, it follows (see for example Lemma 3.1 of \cite{HVM-T-Thas}) that the multiplicity (which might be zero) of the eigenvalue  $\xi_d(\nu_2-r)$ of $PM$ will only depend on $d$, and not on the specific primitive $d$th root of unity. Denote this multiplicity by $a_d$. Analogously the multiplicity (which might be zero) of the eigenvalue  $\xi_d(\nu_3-r)$ of $PM$ will only depend on $d$. Denote this multiplicity by $b_d$. Also, as the sum of all primitive $d$th roots of unity equals $\mu(d)$, where $\mu$ is the M\"obius function, we obtain $$trace(PM)=k-r+ \sum_{d\mid n} a_d\mu(d)(\nu_2-r)+ \sum_{d\mid n} b_d\mu(d)(\nu_3-r).$$ 
On the other hand the trace of $PM$ must equal $-rf+g$, and hence 
\begin{equation}\label{base} k-r+ \sum_{d\mid n} a_d\mu(d)(\nu_2-r)+\sum_{d\mid n} b_d\mu(d)(\nu_3-r)=-rf+g.\end{equation} Setting $r=\nu_3$ we obtain
\begin{equation}\label{eq} k-\nu_3+ \sum_{d\mid n} a_d\mu(d)(\nu_2-\nu_3)=-\nu_3f+g \end{equation}
and 
 \begin{equation}\label{con} k-\nu_3\equiv -\nu_3f+g \pmod{ \sqrt{\Delta}}.\end{equation}
 Finally, in order to compute $a_d+b_d$, we note that as $P$ and $M$ commute and are diagonalizable they are simultaneously diagonalizable (see for example Theorem 1.3.12 in \cite{Matrix}). Let $\mathcal{B}$ be a common eigenbasis for $P$ and $M$, and let $\mathcal{V}_d:=\{v_1,v_2,\hdots,v_m\}\subset\mathcal{B}$ be a basis for the eigenspace of $P$ corresponding to the eigenvalue $\xi_d$ of $P$. If $d=1$ we see that exactly one of the vectors in $\mathcal{V}$ is an eigenvector for the eigenvalue $k-r$ of $M$, exactly $a_d$ of the vectors in $\mathcal{V}$ are  eigenvectors for the eigenvalue $\nu_2-r$ of $M$, and exactly $b_d$ of the vectors in $\mathcal{V}$ are eigenvectors for the eigenvalue $\nu_3-r$ of $M$. Hence $1+a_1+b_1=c$, where $c$ is the number of cycles in the disjoint cycle decomposition of $\phi$. Now let $d\neq1$. Then exactly $a_d$ of the vectors in $\mathcal{V}$ are  eigenvectors for the eigenvalue $\nu_2-r$ of $M$, and exactly $b_d$ of the vectors in $\mathcal{V}$ are eigenvectors for the eigenvalue $\nu_3-r$ of $M$. Hence $a_d+b_d=\sum_{d\mid l} c_l$,  where $c_l$ is the number of cycles of length $l$ in the disjoint cycle decomposition of $\phi$.
 This proves the theorem. \qed
 \medskip
 
This theorem generalizes the previously known Benson type theorems for generalized quadrangles, partial geometries, and partial quadrangles. The key difference with those previous results is that our proof relies on the matrix $M$, which simply is the adjacency matrix plus an arbitrary integer multiple of the identity, rather than relying on the matrix $NN^T$, where $N$ is the incidence matrix of the studied geometry.   This matrix $NN^T$ always equals $A+(t+1)I$, where $A$ is the adjacency matrix of the point graph of the geometry, and $t+1$ is the number of lines through a point in the geometry. But $-t-1$ is not necessarily an eigenvalue of $A$ (for example in the case of partial quadrangles), and hence one cannot guarantee that $NN^T$ has an eigenvalue equal to zero. However, it is exactly the fact that $M$ can be made to have an eigenvalue equal to zero that allows us to deduce Equation (\ref{eq}) and Congruence (\ref{con}) in general. Also the restriction on the value of $a_d+b_d$ is new.
 Though Congruence (\ref{con}) definitely provides a useful tool to analyze possible fixed point structures under automorphisms of a strongly regular graph, we will see in what follows that there are many cases where the more complicated Equality (\ref{eq}) can still be analyzed and  provides  more interesting information. 
 
We also have the following corollary:

\begin{corollary}\label{multiplier}
Let $\cG$ be a strongly regular graph $\srg(v,k,\lambda,\mu)$ with integer eigenvalues, and let $\phi$ be an automorphism of order $n$ of $\cG$. Let $s$ be an integer coprime with $n$. Then $\phi$ and $\phi^s$  map the same number of vertices to adjacent vertices.
\end{corollary}
\proof Because $s$ is an integer coprime with $n$ every vertex fixed by $\phi$ is also fixed by $\phi^s$, and vice versa. Hence $f_\phi=f_{\phi^s}$. Let $P$ be the permutation matrix corresponding to the automorphism $\phi$. As $P$ and $M=A-\nu_3I$ are both diagonalizable and commute, they can be diagonalized simultaneously (see for example Theorem 1.3.12 in \cite{Matrix}). Hence, for every primitive $d$th root of unity $\xi_d$, we can find $a_d^{(\phi)}$ independent vectors $v_1,\hdots,v_{a_d^{(\phi)}}$ that are simultaneously eigenvectors of $M$ with eigenvalue $\nu_2-\nu_3$, and of $P$ with eigenvalue $\xi_d$. Here the superscript $(\phi)$ indicates these are the values $a_d$ corresponding to the automorphism $\phi$. Hence $v_i$, $i=1,\hdots,a_d^{(\phi)}$, is an eigenvector of $P^sM$ with eigenvalue $\xi_d^s(\nu_2-\nu_3)$. As $s$ is coprime to $n$ also $\xi_d^s$ is a primitive $d$th root. Hence $a_d^{(\phi^s)}\geq a_d^{(\phi)}$. But $\phi=(\phi^s)^l$ for some $l$ coprime with $n$, and hence $a_d^{(\phi)}\geq a_d^{(\phi^s)}$. Thus we conclude that $a_d^{(\phi)}=a_d^{(\phi^s)}$. It follows that both $\phi$ and $\phi^s$ produce the same left side in Equation (\ref{eq}). Hence also $g_\phi=g_{\phi^s}$.\qed
\medskip

In Section \ref{diffset} we will see that a well-known multiplier result for partial difference sets is a special case of this corollary.

\subsection{Directed strongly regular graphs}

In this short section we note that both Equalities (\ref{base}), (\ref{eq}) and Congruence (\ref{con}) have a direct analogue for directed strongly regular graphs. We will not explicitly prove this result, but only point out why one can basically copy the proof of Theorem \ref{MAIN} to obtain these analogues.  

A directed strongly regular graph with parameters $(v,k,t,\lambda,\mu)$ is a finite directed graph on $v$ vertices without loops such that every vertex has in-degree and out-degree equal to $k$, each vertex has a constant number $t$ of undirected edges, there are $\lambda$ paths of length 2 between $i$ and $j$ if there is an edge from $i$ to $j$, and there are $\mu$ paths of length $2$ between $i$ and $j$ if there is no edge from $i$ to $j$. Directed strongly regular graphs were introduced by Duval in \cite{Duval}.
The adjacency matrix of a finite directed graph is the square $(0,1)$-matrix $A$ whose columns and rows are labeled by the vertices and is such that $A_{ij}=1$ if and only if there is an edge from vertex $i$ to vertex $j$. The adjacency matrix of a directed strongly regular graph with parameters $(v,k,t,\lambda,\mu)$ satisfies $$A^2+(\mu-\lambda)A-(t-\mu)I=\mu J,$$
$$AJ=JA=kJ.$$ 
Duval (\cite{Duval}) showed that the eigenvalues of the adjacency matrix of a directed strongly regular graph that is not an undirected strongly regular graph or a complete graph are always integers, unless  $A$ is a Hadamard matrix. These eigenvalues are 
$$\kappa_1:=k,$$ 
$$\kappa_2:=\frac{1}{2}(\lambda-\mu+\sqrt{D)}),$$ 
$$\kappa_3:=\frac{1}{2}(\lambda-\mu-\sqrt{D}),$$ 
with respective multiplicities $$m_1:=1,$$ 
$$m_2:=\frac{1}{2}\left(v-1-\frac{2k+(v-1)(\lambda-\mu)}{\sqrt{D}}\right) $$ 
and $$m_3:=\frac{1}{2}\left(v-1+\frac{2k+(v-1)(\lambda-\mu)}{\sqrt{D}}\right),$$
where $D=(\lambda-\mu)^2+4(t-\mu)$.

Let $A$ be the adjacency matrix of a directed strongly regular graph, and let $P$ be the permutation matrix of an automorphism of order $n$ of this graph. As before $PAP^T=A$, and $P^{-1}=P^T$, so that $A$ and $P$ commute. This implies that the eigenvalues of $PM$ are the eigenvalues of $M=A-rI$ multiplied by $n$th roots of unity. Also, if the eigenvalues of $A$ are integers, we are guaranteed that the multiplicity of the eigenvalue $\xi_d \rho$ of $PM$, where $\xi_d$ is a $d$th root of unity and $\rho$ is an eigenvalue of $M$,  will not depend on specific $d$th root of unity. One now easily obtains the following result:

\begin{theorem}
Let $\mathcal{G}$ be a directed strongly regular graph that is not undirected, not a complete graph, and whose adjacency matrix is not a Hadamard matrix.  Let $\phi$ be an automorphism of order $n$ of $\cG$. Let $\mu()$ be the M\"obius function. Then for every integer $r$ and all positive divisors $d$ of $n$, there are non-negative integers  $a_d$ and $b_d$ such that 
\begin{equation}\label{based}k-r +  \sum_{d\mid n} a_d\mu(d)(\kappa_2-r) +  \sum_{d\mid n} b_d\mu(d)(\kappa_3-r)=-rf+g,\end{equation}
where $f$ is the number of fixed vertices of $\phi$ and $g$ the number of vertices that are adjacent to their image under $\phi$.
As a consequence  the following equation and congruence hold:

\begin{equation}\label{eqd}  k-\kappa_3+ \sum_{d\mid n} a_d\mu(d)(\kappa_2-\kappa_3)=-\kappa_3f+g  \end{equation}
and
\begin{equation}\label{cond} k-\kappa_3\equiv -\kappa_3f+g \pmod{ \kappa_2-\kappa_3}. \end{equation}
\end{theorem}
\medskip

Unfortunately there is no nice analogue of Corollary \ref{multiplier}, as one can not guarantee that the adjacency matrix of a directed strongly regular graph is diagonalizable. However, if one knows that $A$ is diagonalizable the result of Corollary \ref{multiplier} remains valid for directed strongly regular graphs.

\section{Automorphisms of prime order}\label{prime}

In this section  we provide the basic idea on how to analyze Equations (\ref{eq}) and (\ref{eqd}). The two examples in Section \ref{examples} will further develop this idea. 
\medskip

As before let $\cG$ be a strongly regular graph $\srg(v,k,\lambda,\mu)$ with integer eigenvalues. Assume $\phi$ is an automorphism of order $p$, prime, of $\cG$. Then we can use the result of theorem \ref{MAIN} to obtain  an interesting divisibility condition. Setting $r=\nu_3$ we obtain $k-\nu_3 +  \sum_{d\mid p} a_d\mu(d)(\nu_2-\nu_3)=-\nu_3f+g.$ As $p$ is prime we have $\mu(1)=1$ and $\mu(p)=-1$. Also $a_1+(p-1)a_p=m_2$. We obtain the following system of linear  equations in $a_1$ and $a_p$:

$$\left\{ \begin{array}{rcl}  a_1(\nu_2-\nu_3)-a_p(\nu_2-\nu_3) & = & -\nu_3(f-1)+g-k \nonumber \\ a_1+(p-1)a_p & = & m_2. \nonumber \end{array} \right. $$
This system can easily be solved for $a_1$ and $a_p$  
which should be non-negative integers.  This can be of particular interest if one is trying to disprove the existence of a hypothetical strongly regular graph which should admit certain automorphisms, for example a certain Cayley graph (see Sections \ref{diffset} and \ref{examples}). In what follows we will see how in certain cases it is possible to obtain results that go beyond automorphisms of prime order.
\medskip

A similar result can be obtained for directed strongly regular graphs in the obvious way. 

\section{Partial difference sets in abelian groups}\label{diffset}

We will now focus on strongly regular graphs where the existence of a certain automorphism is a priori known. Strongly regular Cayley graphs on abelian groups are important graphs is this category. In this case the graph is equivalent to a so-called partial difference set (PDS). We will recover some known results on these objects with new elementary proofs, as well as some new results. We first review some basic definitions.
\medskip

Let $G$ be a finite abelian group of order $v$, and let $\cD$ be a $(v,k,\lambda,\mu)$ partial difference set in $G$, that is, $\cD$ is a $k$-subset of $G$ with the property that the expressions $\phi\psi^{-1}$, $\phi,\psi\in \cD$, represent each nonidentity element in $\cD$ exactly $\lambda$ times, and each nonidentity element of $G$ not in $\cD$ exactly $\mu$ times. Further assume that $\cD^{(-1)}=\cD$ (that is, if $\phi\in \cD$, then so is $\phi^{-1}$) and $e\notin\cD$, where $e$ is the identity in $G$, that is, $\cD$ is a so-called regular partial difference set. Then it is well known (see e.g. \cite{MA94}) that the Cayley graph $\cG:=(G,\cD)$ is a strongly regular graph with parameters srg$(v,k,\lambda,\mu)$. Recall that two elements $\phi$ and $\psi$ of $G$ are adjacent in $\cG$ if and only if $\psi\phi^{-1}\in\cD$. Also, if $\lambda\neq\mu$ then $\cD^{(-1)}=\cD$ is automatically fulfilled (see \cite{MA94}).
\medskip

Let $\cG=(G,\cD)$ be a strongly regular Cayley graph on an abelian group $G$, and assume that $\cG$ is not a conference graph. Then the nonidentity elements of $G$ act in an obvious way as fixed point free automorphisms on $\cG$. Furthermore, such automorphism, $\phi$, maps either every vertex to an adjacent vertex (if and only if $\phi\in \cD$), or maps no vertex to an adjacent vertex. The latter can be seen as follows. Assume $\phi$ maps vertex $\psi$ to adjacent vertex $\psi\phi$. Then clearly $(\psi\phi)\psi^{-1}=\phi\in \cD$. For every other vertex $\gamma$ we then have that $(\gamma\phi)\gamma^{-1}=\phi\in\cD$, and hence every vertex is mapped to an adjacent vertex. We obtain the following result:

\begin{corollary}
With the above notation, it holds that $$2k-\lambda+\mu\equiv0\equiv v \pmod{ \sqrt{\Delta}}.$$
\end{corollary}
\proof First let $\phi$ be an element of $G\setminus(\cD\cup\{e\})$. Then we can apply Congruence \ref{con} with $f=0$ and $g=0$. We obtain \begin{equation}\label{help}k-\frac{1}{2}(\lambda-\mu-\sqrt{\Delta})\equiv0\pmod{\sqrt{\Delta}}.\end{equation}
Next, let $\phi$ be an element of $\cD$. Then we can apply Congruence \ref{con} with $f=0$ and $g=v$. We obtain $$k-\frac{1}{2}(\lambda-\mu-\sqrt{\Delta})\equiv v\pmod{\sqrt{\Delta}}.$$ Combining both Congruences proves the result. \qed
\medskip

This congruence was obtained with a different proof by S.L. Ma \cite{MA94}.

Corollary \ref{multiplier} can be used to obtain a multiplier result that was previously obtained by Ma \cite{MA94}. 

\begin{corollary}\label{multi}
Let $\cD$ be a regular $(v,k,\lambda,\mu)$ partial difference set not containing the identity in the abelian group $G$. Furthermore assume $\Delta$ is a perfect square. Then $\cD^{(s)}=\cD$ for every $s$ coprime to $v$.
\end{corollary}
\proof This follows immediately from Corollary \ref{multiplier} and the fact that an element of $G$ belongs to $\cD$ if and only if it maps every vertex to an adjacent vertex in the corresponding Cayley graph.\qed

\medskip

We next look at Equality (\ref{eq}). It is obvious that one obtains two equalities, one for $\phi\notin \cD$, and one for $\phi\in\cD$.  

If we assume  that $\phi$ has prime order the results from Section \ref{prime} can be applied, and we obtain:

If $\phi\notin\cD$ then $a_p=\frac{-v(\lambda-\mu-\sqrt{\Delta)}}{2p\sqrt{\Delta}}$ must be an integer. If $\phi\in\cD$ then $a_p=\frac{-2v-v(\lambda-\mu-\sqrt{\Delta)}}{2p\sqrt{\Delta}}$ must be an integer. At first sight this does not seem to provide much information (because $p$ definitely divides $v$), however, if we use this as the basis for studying automorphism of higher order useful conclusions can be drawn (see Section \ref{examples}).

\section{Two examples}\label{examples}

We will illustrate how the ideas developed in this paper can be used by disproving the existence of a $(100,33, 8, 12)$ and a $(100,36, 14, 12)$ regular partial difference set in an Abelian group of order $100$. These were the only two parameter sets for regular PDS in Abelian groups of size at most $100$ for which (non)-existence had not been settled (see \cite{MA97}, \cite{MA-MOK}).  It is worthwhile to note that a  $(100,36, 14, 12)$ PDS does exist in non-Abelian groups of order $100$ (see \cite{KLin}).
We start with the $(100,33, 8, 12)$ case.
\medskip

Assume by way of contradiction that $\cD$ is a $(100,33, 8, 12)$ regular PDS in the Abelian group $G$. We first show that $G$ cannot contain a subgroup of order $25$ by analyzing elements of order $5^k$.

{\sc Elements of order $5$.} Let $\phi$ be an element of order $5$ in $G$. We compute $a_1^{(\phi)}$ and $a_5^{(\phi)}$. 

$$\begin{array}{|c|c|c|} \hline o(\phi)=5 &  a_1^{(\phi)} & a_5^{(\phi)} \\ \hline \phi\in\cD  & 18&12 \\ \hline \phi\notin\cD & 10&14 \\ \hline \end{array}$$

{\sc Elements of order $25$.} Next assume $\phi$ is an element of order $25$ in $G$. Then $\phi^5$ is an element of order $5$, and we obtain (using the same ideas as in the proof of Corollary \ref{multiplier}):

$$a_1^{(\phi^5)}=a_1^{(\phi)}+4a_5^{(\phi)}, \ \ \mathrm{ and }\ \  4a_5^{(\phi^5)}=20a_{25}^{(\phi)}.$$

 From the last equation we see that $a_5^{(\phi^5)}$ must be a multiple of $5$, however, a quick inspection of the table with multiplicities for elements of order $5$ shows this is never the case. Hence $G$ does not contain elements of order $25$, and $G\cong\mathbb{Z}_2\times\mathbb{Z}_2\times\mathbb{Z}_5\times\mathbb{Z}_5$ or $G\cong\mathbb{Z}_4\times\mathbb{Z}_5\times\mathbb{Z}_5$.

\medskip

By Corollary \ref{multiplier} the elements of order $5$, $10$, and $20$ in $\cD$ come in sets of $4$, $4$  and $8$ respectively. Elements of order $4$ in $\cD$ come in sets of two. As the size of $\cD$ is $33$ it now easily follows that $\cD$ contains a unique element of order $2$, while all other elements in $\cD$ have order divisible by $5$. We now look at the elements of $G$ whose order is divisible by $2$.

{\sc Elements of order $2$} If $\phi$ is an element of order $2$ in $G$ then we obtain the following multiplicities:

$$\begin{array}{|c|c|c|} \hline o(\phi)=2 & a_1^{(\phi)} & a_2^{(\phi)}  \\ \hline \phi\in D &36 &30 \\ \hline \phi\notin D &31 &35 \\ \hline  \end{array}$$

{\sc Elements of order $10$} Now let $\phi$ be an element in $G$ of order $10$. Using that $\phi^5$ has order $2$ and $\phi^2$ order $5$ we obtain

$$a_1^{(\phi^5)}=a_1^{(\phi)}+4a_5^{(\phi)}, \ \ \mathrm{ and }\ \  a_2^{(\phi^5)}=a_2^{(\phi)}+4a_{10}^{(\phi)},$$

$$a_1^{(\phi^2)}=a_1^{(\phi)}+a_2^{(\phi)}, \ \ \mathrm{ and }\ \  4a_5^{(\phi^2)}=4a_5^{(\phi)}+4a_{10}^{(\phi)},$$

and finally from Equation (\ref{eq}) 

$$ 40 +(a_1^{(\phi)}-a_2^{(\phi)}-a_5^{(\phi)}+a_{10}^{(\phi)})10=g,$$
where $g=0$ or $g=100$ depending on whether $\phi\not\in\cD$ or $\phi\in\cD$. There are eight cases to consider depending on whether $\phi,\phi^2, \phi^5$ all or not belong to $\cD$.

We obtain 

$$\begin{array}{|c|c|c|c|c|} \hline o(\phi)=10& a_1^{(\phi)} & a_2^{(\phi)} & a_5^{(\phi)} & a_{10}^{(\phi)} \\ \hline \phi^5\in\cD, \phi^2\in\cD,\phi\in\cD &12 &6 & 6&6 \\ \hline \phi^5\in\cD, \phi^2\in\cD,\phi\notin\cD &8 &10 & 7&5 \\ \hline \phi^5\in\cD, \phi^2\notin\cD,\phi\in\cD &8 &2 & 7&7 \\ \hline \phi^5\in\cD, \phi^2\notin\cD,\phi\notin\cD &4 &6 & 8&6 \\ \hline \phi^5\notin\cD, \phi^2\in\cD,\phi\in\cD &11 &7 & 5&7 \\ \hline \phi^5\notin\cD, \phi^2\in\cD,\phi\notin\cD &7 &11 & 6&6 \\ \hline \phi^5\notin\cD, \phi^2\notin\cD,\phi\in\cD &7 &3 & 6&8 \\ \hline \phi^5\notin\cD, \phi^2\notin\cD,\phi\notin\cD &3 &7 & 7&7 \\ \hline

\end{array}$$

First note that from this table we see that possibilities 1, 5 and 6 cannot occur. The reason for this is that $\phi$ acts with $10$ orbits of length $10$, and hence has each $10$th root of unity as an eigenvalue with multiplicity exactly $10$. This implies that none of the $a_i^{(\phi)}$, $i=1,2,5,10$, can be larger then $10$. 

We now exclude $G\cong\mathbb{Z}_4\times\mathbb{Z}_5\times\mathbb{Z}_5$  and $G\cong\mathbb{Z}_2\times\mathbb{Z}_2\times\mathbb{Z}_5\times\mathbb{Z}_5$ one by one. First assume $G\cong\mathbb{Z}_4\times\mathbb{Z}_5\times\mathbb{Z}_5$. Then $G$ would contain an element $\phi$ of order $20$, yielding 

$$4a_{10}^{(\phi^2)}=8a_{20}^{(\phi)}, \ \ \mathrm{ andÊ}\ \  a_{2}^{(\phi^2)}=2a_{4}^{(\phi)}.$$ As every element of order $10$ is the square of an element of order $20$  this further excludes, by the obvious parity argument, possibilities 2, 3, 7 and 8 from the above table. Hence only possibility 4 is left. As every element of order $5$ is the square of an element of order $10$, it follows that $\cD$ does not contain any element of order $5$ or $10$, and hence $\cD$ consists of the unique element of order $2$ and $32$ elements of order $20$. However, it is easy to check that the elements of order $4$ can never be written as a difference of two elements of order $20$ or as the difference of an element of order $20$ and the unique element of order $2$. This contradicts $\cD$ is a partial difference set with the given parameters. 

Finally assume that $G\cong\mathbb{Z}_2\times\mathbb{Z}_2\times\mathbb{Z}_5\times\mathbb{Z}_5$.
We already know that $\cD$ must contain a unique element of order $2$, say $\iota$. Now assume that $\gamma$ would be an element of order $5$ in $\cD$, and let $\kappa$ be an element of order $2$ not in $\cD$. Let $\gamma'$ be the unique element such that $\gamma'^2=\gamma$. Then $\phi=\kappa\gamma'$ has order $10$ and $\phi^2=\gamma\in \cD$, whereas $\phi^5=\kappa\notin\cD$. Hence we must be in case 5 or 6 of our table. However, these cases were previously excluded. This implies that $\cD$ cannot contain any element of order $5$. For convenience let us write the elements of $G$ as $(a,b,c,d)$, where addition is done component wise modulo $2$ in the first two components and modulo $5$ in the two last components. Without loss of generality we can assume that $\iota=(1,1,0,0)$. All other elements of $\cD$ are of the form $(1,0,a,b), (1,1,a,b)$ or $(0,1,a,b)$, where $(a,b)\neq(0,0)$. It is easy to see that the only way to write $(1,0,0,0)$ as a difference of elements in $\cD$ is as $(1,0,0,0)=(1,1,a,b)-(0,1,a,b)=(0,1,a,b)-(1,1,a,b)$. However, as the elements of order $10$ in $\cD$ come in groups of $4$, this implies that the differences in $\cD$ that produce $(1,0,0,0)$ come in sets of $8$. But the $\mu$-parameter for this partial difference set equals $12$, not a multiple of $8$, the final contradiction.

\medskip

Next we exclude the existence of a $(100,36, 14, 12)$ regular PDS in an Abelian group $G$. As the ideas and techniques are very similar we will be briefer  in our arguments.

{\sc Elements of order $5$.}

$$\begin{array}{|c|c|c|} \hline o(\phi)=5 &  a_1^{(\phi)} & a_5^{(\phi)} \\ \hline \phi\in\cD  & 12&6 \\ \hline \phi\notin\cD & 4&8 \\ \hline \end{array}$$.

Because $a_5^{(\phi)}$ is never divisible by $5$ this again excludes the existence of a subgroup of order $25$ in $G$. 

It is useful to make an argument based on Corollary \ref{multiplier} at this point. As there are only $3$ elements in $G$ whose order is a power of $2$, the elements of order divisible by $5$ in $\cD$ come in sets whose size is a multiple of $4$, and $\left|\cD\right|$ is divisible by $4$, it follows that $\cD$ does not contain any element of order a power of $2$. This allows us to only compute partial tables of multiplicities.

{\sc Elements of order $2$, $4$ and $10$.}

$$\begin{array}{|c|c|c|} \hline o(\phi)=2 & a_1^{(\phi)} & a_2^{(\phi)}   \\ \hline \phi\notin D &16 &20 \\ \hline  \end{array}$$

$$\begin{array}{|c|c|c|c|} \hline o(\phi)=4 & a_1^{(\phi)} & a_2^{(\phi)} & a_4^{(\phi)} \\ \hline \phi\not\in\cD &6 &10 & 10 \\ \hline \end{array}$$

$$\begin{array}{|c|c|c|c|c|} \hline o(\phi)=10& a_1^{(\phi)} & a_2^{(\phi)} & a_5^{(\phi)} & a_{10}^{(\phi)}    \\ \hline \phi^5\notin\cD, \phi^2\in\cD,\phi\in\cD &8 &4 & 2&4 \\ \hline \phi^5\notin\cD, \phi^2\in\cD,\phi\notin\cD &4 &8 & 3&3 \\ \hline \phi^5\notin\cD, \phi^2\notin\cD,\phi\in\cD &4 &0 & 3&5 \\ \hline \phi^5\notin\cD, \phi^2\notin\cD,\phi\notin\cD &0 &4 & 4&4 \\ \hline

\end{array}$$

Now first assume $G\cong\mathbb{Z}_4\times\mathbb{Z}_5\times\mathbb{Z}_5$, and let $\phi$ be an element of order $20$. Then, as before, we see that $a_2^{(\phi^2)}$ and $a_{10}^{(\phi^2)}$ must be even. Hence only case 1 ($\phi^2\in\cD$) and  case 4 ($\phi^2\notin\cD$) in the previous table can occur. We obtain

$$\begin{array}{|c|c|c|c|c|c|c|} \hline o(\phi)=20& a_1^{(\phi)} & a_2^{(\phi)} & a_4^{(\phi)} & a_5^{(\phi)} & a_{10}^{(\phi)}  & a_{20}^{(\phi)}  \\ \hline \phi^2\in\cD, \phi\in\cD &6 &2 & 2&0&2&2 \\ \hline \phi^2\in\cD, \phi\notin\cD &2 &6 & 2&1&1&2 \\ \hline \phi^2\notin\cD, \phi\in\cD &2 &-2 & 2&1&3&2 \\ \hline , \phi^2\notin\cD,\phi\notin\cD &-2 &2 & 2&2&2&2 \\ \hline
\end{array}$$

As $\phi$ acts with $20$ orbits of length $5$, it has each $20$th root of unity as an eigenvalue with multiplicity exactly $5$. This implies that none of the $a_i^{(\phi)}$, $i=1,2,4,5,10, 20$, can be larger then $5$.  This excludes the first two possibilities. Finally the last two possibilities are excluded as multiplicities have to be non-negative integers.

Finally assume that $G\cong\mathbb{Z}_2\times\mathbb{Z}_2\times\mathbb{Z}_5\times\mathbb{Z}_5$.
Let $\iota$ be any element of order $2$. If $\iota=\alpha\beta^{-1}$, for $\alpha,\beta\in\cD$, then also $\iota=\beta\alpha^{-1}$. As the elements of $\cD$ have order $5$ or $10$ they come in sets of $4$. Hence the number of ways in which $\iota$ can be written as a difference of elements of $\cD$ is a multiple of $8$. However, $\mu=12$, the final contradiction.

\section{Projective two-weight sets}

In this last section we will shortly describe how the results of this paper can be used in the analysis of hypothetical  projective two-weight set with given parameters.

Let PG$(n,q)$ be the $n$-dimensional projective space over the field $\mathbb{F}_q$. A projective two-weight set $K$ of type $(w_1,w_2)$, $ w_1\neq w_2\neq0$, in PG$(n,q)$ is a set of points in PG$(n,q)$ with the property that every hyperplane of PG$(n,q)$ intersects $K$ in either $w_1$ or $w_2$ points. It is well known that such sets are equivalent to linear two-weight codes and give rise to strongly regular graphs (see for example \cite{CalKant}). Here we will only describe how these sets give rise to a strongly regular graph, and then apply the results earlier obtained in the paper.

From here on, let $K$ be a two-weight set of type $(w_1,w_2)$ and size $N$ in PG$(n,q)$. Embed PG$(n,q)$ as a hyperplane in PG$(n+1,q)$. Construct the graph $\cG$ as follows: the vertices of $\cG$ are the points of PG$(n+1,q)\setminus$PG$(n,q)$; two vertices $g$ and $h$  of $\cG$ are adjacent if and only if the projective line $<g,h>$ intersects PG$(n,q)$ in a point of $K$. Then it is well known that $\cG$ is strongly regular with parameters $v=q^{n+1}$, $k=(q-1)N$, $\lambda=k^2+3k-(k+1)q(2n-w_1-w_2)$, $\mu=(n-w_1)(n-w_2)/q^{n-1}$ (see \cite{CalKant}). It is obvious that the group $T$ of translations of PG$(n+1,q)$ with axis PG$(n,q)$ acts as a sharply transitive abelian group of automorphisms.  Now fix any point (vertex of $\cG$) $g$ of PG$(n+1,q)\setminus$PG$(n,q)$, and identify any point $h$ of PG$(n+1,q)\setminus$PG$(n,q)$ with the unique element $\tau$ of $T$ with the property that $g^\tau=h$. If under this correspondence $\cD$ is the subset of $T$ corresponding to the vertices of $\cG$ adjacent to $p$, then it is clear that $\cD$ is a regular PDS in $T$, and that $\cG=(T,\cD)$. Hence all results from Section \ref{diffset} can be applied. However, the graph $\cG$ admits another interesting group of automorphisms. Again, let $g$ be a fixed chosen point of PG$(n+1,q)\setminus$PG$(n,q)$. Then the group $H$ of homologies with center $g$ and axis PG$(n,q)$ will act obviously as a group of automorphisms of $\cG$.  Any non-identity element of $H$ will have exactly one fixed vertex ($g$) and map $k=(q-1)N$ vertices to adjacent vertices.  First note that simply applying Congruence (\ref{con}) does not produce any useful results (it simply yields $k-\nu_3\equiv k-\nu_3\pmod{\nu_2-\nu_3}$). However, just as was discussed in Sections \ref{diffset} and \ref{examples}, analyzing Equation (\ref{eq}) does provide new information. The group $H$ has order $q-1$, and for every prime divisor $p$ of $q-1$ we can pick an element $\phi$ of $H$ of that order and apply the result of Section \ref{prime}. We obtain that $a_1^{(\phi)}=a_p^{(\phi)}$ and $a_1^{(\phi)}=m_2/p$. As $p$ divides $v-1=q^{n+1}-1$ this typically is not a strong condition. However, it might, just as in Section \ref{examples}, be used as the basis to analyze automorphisms of higher order in $H$.

\subsubsection*{Acknowledgement} The authors want to thank Tim Penttila for suggesting to look at directed strongly regular graphs.

\end{document}